\documentclass[11pt]{amsart}
 
\usepackage{amsmath}
\usepackage{amsfonts}
\usepackage{amssymb}
\usepackage{epsfig}
\usepackage{epstopdf}
\usepackage{color}

  \newtheorem{theorem}{Theorem}
  \newtheorem{lemma}[theorem]{Lemma}
  
  \newtheorem{problem}[theorem]{Problem}

  \newtheorem{proposition}[theorem]{Proposition}

\newtheorem{conjecture}[theorem]{Conjecture}

  \theoremstyle{definition}
  
  \newtheorem{example}[theorem]{Example}

\newcommand{\CC}{\mathbb C}
\newcommand{\RR}{\mathbb R}

\newcommand{\PP}{\mathbb P}

\newcommand{\TP}{\mathbb{TP}}
\newcommand{\TT}{\mathbb T}
 
\newcommand{\C}{\mathcal{C}}
 \newcommand{\cT}{\mathcal{T}}
 \newcommand{\B}{\mathcal{B}}

\title{Representing Tropical Linear Spaces by Circuits}



\author{Josephine Yu}
\address{Department of Mathematics, University of California, Berkeley, CA 94720.}
\email{jyu@math.berkeley.edu}

\author{Debbie S. Yuster}
\address{Department of Mathematics, Columbia University,
2990 Broadway,
509 Mathematics Building,
Mail Code: 4406,
New York, NY 10027 
}
\email{debbie@math.columbia.edu}

\date{November 19, 2006}

\begin{document}

\begin{abstract}
We study representations of tropical linear spaces as intersections of tropical hyperplanes of circuits.  For several classes of matroids, we describe minimal tropical bases.  We also show that every realizable tropical linear space has a natural, tropically linear parametrization involving its cocircuits.
\end{abstract}

 \maketitle
 
 \section{Introduction}

 Let $\TT = (\RR \cup \{\infty\}, \oplus, \odot )$ be the \emph{tropical semiring} where $\oplus$ is taking minimum and $\odot$ is the usual addition.    We can mod out $\TT^n$ by tropical scalar multiplication to get the \emph{tropical projective space} $\TP^{n-1} = (\TT^n \backslash (\infty, \dots, \infty) ) / \RR(1,\dots, 1)$, which is sometimes more convenient to work in.  

 Tropical linear spaces are tropical analogues of usual linear spaces.  The \emph{tropical hyperplane} of a tropical linear form $(c_1 \odot x_1) \oplus \cdots \oplus (c_n \odot x_n), c \in \TP^{n-1},$ is the set of points $x \in \TP^{n-1}$ such that the minimum in the linear form is attained at least twice.  Let us recall the setup from \cite{Sp}.  Let $d \leq n$ be positive integers.  A point $p \in \TP^{{[n] \choose d}-1}$ is called a \emph{tropical Pl\"{u}cker vector} if for every $d-2$-subset $S$ of $[n]$ and four distinct elements $i,j,k,l \in [n] \backslash S$, the minimum in $(p_{Sij} \odot p_{Skl}) \oplus  (p_{Sik} \odot p_{Sjl}) \oplus(p_{Sil} \odot p_{Sjk})$ is attained at least twice.  
 Given a tropical Pl\"{u}cker vector $p \in \TP^{{[n] \choose d}-1}$,  for each $d+1$-subset $I \subset [n]$ we can define a tropical linear form $\oplus_{i \in I}  \left( p_{I\backslash \{i\}} \odot x_i \right)$ called a \emph{circuit}.    The \emph{tropical linear space} corresponding to $p$ is the intersection of the tropical hyperplanes of these circuits.  A \emph{tropical basis} of a tropical linear space is a set of defining linear forms for the space. Tropical bases are not unique, and need not be minimal in any sense. Much of this paper is concerned with finding \emph{minimal} tropical bases. Tropical linear spaces are needed to compute tropical discriminants \cite{DFS}, and thus a minimal basis is desirable.

A tropical linear space whose defining tropical linear forms have coefficients all 0 or $\infty$ is called \emph{constant coefficient}. Its associated tropical hyperplanes are determined by the supports of the tropical linear forms, that is, the entries with non-$\infty$ coefficients. As a result, the conditions for being a tropical basis depend only on those supports.  
In this case, we can deal with the matroid whose dependent sets are supports of these tropical linear forms.  In Sections 2 and 3, we deal exclusively with the constant coefficient case, which amounts to finding tropical bases of matroids. We describe minimal tropical bases for several classes of matroids. Our main findings, appearing in Section 3, are that graphic matroids, cographic matorids, and the matroid $R_{10}$ have \emph{unique} minimal tropical bases. We hope to extend these findings to all regular matroids.

In Section 4, our main result concerns the tropical rank of a matrix whose rows form a tropical basis. Furthermore, we conjecture a criterion for being a tropical basis in the non-constant case.  Finally, in Section \ref{sec:param}, we show that there is a natural parametrization of tropical linear spaces in terms of cocircuits.

\section{The constant coefficient case}

As discussed in the Introduction, constant coefficient tropical linear spaces can be described in terms of their associated matroids. We will deal with the matroid whose circuits are supports of circuits of a tropical linear space.  

Let $M$ be a matroid and $\C$ be the collection of its circuits.
For a circuit $C \in \C$, let $\cT(C)$ be the set of points $x \in \TP^{n-1}$ such that the minimum value in $\{x_i : i \in C\}$ is attained at least twice.   The set $\cT(\C) := \bigcap_{C \in \C} \cT(C)$ is a polyhedral fan called the {\em tropical variety} or the {\em Bergman fan}  of $M$.  Given a subset $\B\subset\C$, define $\cT(\B):=\bigcap_{C\in\B} \cT(C)$. The set $\B$ is called a {\em tropical basis} of $M$ if $\cT(\B) = \cT(\C)$.

\begin{problem} 
Identify a minimal tropical basis for any matroid.
\end{problem}

It was shown in \cite{AC} that the intersection of the tropical variety of a matroid $M$ and a sphere centered at the origin is a geometric realization of the order complex of the lattice of flats of $M$.  In a matroid, an element $e$ is called a \emph{loop} if $\{e\}$ is a circuit, and two elements $e_1, e_2$ are said to be \emph{parallel} if $\{e_1, e_2\}$ is a circuit.
 Since removing the loops and replacing each parallel class with a single element in a matroid does not change the lattice of flats, we may assume that our matroids are simple, i.e.\ contain no loops or parallel elements.  Since the circuits of the direct sum (or 1-sum) of two matroids is the union of circuits of the summands, the following is clear.

\begin{lemma}\label{directsum}
If matroids $M$, $M_1$, and $M_2$ are such that $M = M_1 \oplus M_2$, then the tropical bases of $M$ are precisely unions of tropical bases of $M_1$ and $M_2$.
\end{lemma}

Since every matroid is the direct sum of its connected components, we can restrict attention to connected matroids.

Each circuit in the tropical basis ``excludes'' certain points from being in the tropical variety, namely those values which induce a unique minimum on the terms of that circuit. We can cut down a tropical basis of a matroid as long as the smaller circuit set excludes the same points as the larger one. In comparing the excluded points, it suffices to consider the 0/1 points:

\begin{lemma}
For any $\B \subset \C$,
$$[\cT(\B) \cap \{0,1\}^n = \cT(\C) \cap \{0,1\}^n] \implies [\cT(\B) = \cT(\C)].$$ In other words, to test equality of $\cT(\B)$ and $\cT(\C)$, it is sufficient to check that they agree on 0/1 points.
\end{lemma}

\begin{proof}
Suppose we are given a subset $\B\subset\C$ of circuits of a matroid $M$, and we would like to know if $\B$ is a tropical basis of $M$. Since each circuit of a matroid excludes points, $\cT(\B)$ contains, or is equal to, $\cT(\C)$. Thus, we need only worry about ``extra'' points, i.e.\ points which are in $\cT(\B)$ but not in $\cT(\C)$. We show that if $\cT(\B)$ contains an extra point $x$, then $\cT(\B)$ contains an extra 0/1 point. Thus checking that there are no extra 0/1 points shows equality of $\cT(\B)$ and $\cT(\C)$.

Suppose $\cT(\B)$ contains an extra point $x$. Then there exists a circuit $C\in \C\backslash\B$ which excludes $x$, and therefore the set $\{x_i:i\in C\}$ has a unique minimum $m$. From $x$, construct the 0/1 vector $v$ as follows: 
$$
v_i = \left\{
\begin{array}{cc}
1 & \text{ if } x_i > m \\
0 & \text{ if } x_i \leq m
\end{array} 
\right.
$$
We see that $v$ is excluded by $C$, so $v\notin\cT(\C)$. However, for any circuit $\tilde C$ such that $\{x_i:i\in \tilde C\}$ attains its minimum at least twice, $\{v_i:i\in \tilde C\}$ also attains its minimum at least twice, so $v\in\cT(\B)$. 
\end{proof}

So it is possible to remove some circuits from a tropical basis of a matroid, if the same 0/1 vectors are excluded in the smaller set.   It was shown in \cite{AC} that a 0/1-vector is in the tropical variety if and only if its support, the set of coordinates with non-zero values, is a flat of the matroid.  Hence the excluded points are non-flats.  A collection of circuits is a tropical basis if and only if it excludes all 0/1 non-flats.
\\

If two circuits $C_1, C_2$ have a unique element in their intersection, \emph{pasting} them means taking their symmetric difference $C_1 \triangle C_2 = (C_1 \backslash C_2) \cup (C_2 \backslash C_1)$.

\begin{lemma}\label{pasting}
If a collection $\mathcal{S}$ of circuits of a matroid has the property that every other circuit of the matroid can be obtained by successively pasting circuits in $\mathcal{S}$, then $\mathcal{S}$ is a tropical basis.
\end{lemma}

\begin{proof} 

Suppose we have two circuits, $A$ and $B$, and a weight assignment such that $A$ and $B$ both have their minimum value attained at least twice. Form the circuit $C$ by pasting $A$ and $B$ together along some element $e$. We claim that $C$ too attains its minimum at least twice.

The minima may be attained in one of three ways:

Case 1) $A$ attains its minimum on $a$ and $a'$, neither of which is $e$, and $B$ attains its minimum on $b$ and $b'$, neither of which is $e$. In this case, $C$ attains its minimum twice, on either $a$ and $a'$, or $b$ and $b'$ (or all four), depending on which edge weights are minimal.

Case 2) $A$ attains its minimum on edges $a$ and $e$ and $B$ attains its minimum on edges $b$ and $b'$, neither of which is $e$. Then $b$ and $b'$ have lower weight than $e$, so $C$ attains its minimum on $b$ and $b'$.

Case 3) A attains its minimum on edges $a$ and $e$, and $B$ attains its minimum on edges $b$ and $e$. Then $weight(a)=weight(e)=weight(b)$, and thus $C$ attains its minimum on $a$ and $b$.
\end{proof}

\subsection{Partition Matroids}

Given a partition $\Pi$ of $[n]$, the \begin{it}partition matroid\end{it} $M_\Pi$ is the matroid whose circuits are pairs of elements in the same block of the partition.  This is the case when the defining ideal of the linear space is also binomial.

The minimal tropical basis consists of enough pairs in each block to form a ``spanning tree'' on that block. This ensures uniform weighting within each block, and thus forces the minimum to be attained twice on each circuit.

\subsection{Uniform Matroids}

The uniform matroid $U_{d,n}$ is the matroid arising from a generic set of $n$ points in $\mathbb{R}^d$. The circuits of $U_{d,n}$ are the $(d+1)$-subsets of $[n]$.  Any tropical basis of $U_{d,n}$ must contain at least $\frac{1}{d+1} {n \choose d}$ elements \cite[Theorem 2.10]{BJSST}, and the bound is not tight.

\begin{lemma}\label{uniform}
An inclusion minimal tropical basis $\B$ for a uniform matroid $U$ is given by
$$ \B=\{C:i\in C\},$$ where $i$ is any fixed element of the matroid.  
\end{lemma}

\begin{proof}
We must show that given a point $x\in\mathbb{R}^n$ that is excluded by some circuit of $U$, we can find a circuit in $\B$ which excludes $x$. As shown above, we can restrict our attention to 0/1 points. Without loss of generality let us fix $i=1$, so that the tropical basis $\B$ consists of all circuits of $U$ containing the element $1$. Consider a 0/1 point $x\in\mathbb{R}^n$ that is excluded by some circuit $C_x$. The minimum of $\{x_j:j\in C_x\}$ is attained uniquely.  Consequently, there is exactly one element $k$ of $C_x$ such that $x_k=0$.

If $C_x$ contains the element $1$, then it is in $\B$ and we're done. Otherwise, consider the following two cases:

Case 1 ($x_1=1$): Let $\tilde{C_x}$ be the circuit obtained from $C_x$ by replacing any element of $C_x$ \begin{it}other than\end{it} $k$ by the element $1$. 

Case 2 ($x_1=0$): Let $\tilde{C_x}$ be the circuit obtained from $C_x$ by replacing $k$ with $1$.

In both cases,  $\tilde{C_x}$ is a circuit of $U$, since every $(d+1)$-subset of $[n]$ is a circuit. It contains the element 1, so it is in $\B$. Finally, $\tilde{C_x}$ excludes $x$, since the minimum of $\{x_j:j\in \tilde{C_x}\}$ is attained uniquely. Thus $x$ is excluded by the tropical basis $\B$.

To see that $\B$ is inclusion minimal, consider a circuit $C^*\in\B$. The point $x$, where $x_i=1$ for all $i\in C^*\backslash 1$ and $x_i=0$ otherwise, would be excluded by $\B$, but not by $\B\backslash \{C^*\}$. Thus it is necessary to include $C^*$ in this tropical basis.
\end{proof}

Note that the tropical bases given for uniform matroids are not unique. For example, consider $U_{2,4}$. It has four circuits, namely 123, 124, 134, and 234. Any collection of three out of the four circuits forms an inclusion- and cardinality- minimal tropical basis. Yet, there is no single circuit that \begin{it}must\end{it} be in a tropical basis of $U_{2,4}$. Additionally, the tropical bases given above for uniform matroids are not in general cardinality-minimal, despite being inclusion-minimal. For example, the smallest tropical bases of $U_{2,5}$ contain five elements (123, 124, 125, 134, and 345 for example), however there are inclusion-minimal tropical bases containing six elements, namely the six circuits which contain the element `1'.

\subsection{Fano plane}

\begin{figure}[htb]
\includegraphics[scale=0.6]{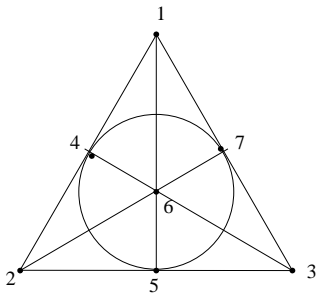}
\caption{Fano plane}
\label{fanoplane}
\end{figure}

The Fano plane gives rise to a matroid whose ground set consists of its 7 vertices. The circuits are the ``lines'' in the point configuration (see Figure \ref{fanoplane}) and any 4-subset not containing a line.

The unique minimal tropical basis for the Fano plane consists of the seven 3-element circuits, one arising from each line (124, 137, 156, 235, 267, 346, \& 457).

The lines are necessary: If one of the lines, $L$, were to be excluded, we could assign weight 1 to two of the line's points, weight 0 to the remaining point, and weight 0 to all points off the line. Thus every circuit except $L$ would have its minimum attained twice, while $L$ would have a unique minimum.

The lines are sufficient: Each 4-element circuit (1236, 1467, 2456, 3567, 1257, 1345, and 2347) can be obtained by pasting two 3-element circuits together along a shared element, and deleting that element. For example, we can think of the circuit 1236 as being pasted from the circuits 124 and 346. Thus, by Lemma \ref{pasting}, the 3-element circuits are sufficient.

\begin{problem}
Find a minimal tropical basis for transversal matroids \cite[Section 3.3]{STY}.
\end{problem}

\begin{problem}
What happens to tropical bases under taking minors (deletion and contraction)?
\end{problem}

\section{Regular matroids}

A \emph{regular matroid} is one that is representable over every field. In this section we characterize minimal tropical bases of certain important classes of regular matroids, namely graphic matroids, cographic matroids, and the matroid $R_{10}$.  Seymour showed that every regular matroid can be constructed by piecing together matroids of these types \cite{seymour}.

\subsection{Graphic Matroids}

A graphic matroid is formed from a graph $G$. The edges of $G$ form the ground set, and the circuits of the matroid are the edge collections corresponding to cycles of $G$, where a cycle is a closed walk all of whose vertices have degree two in the cycle. If $G$ contains $n$ edges, we can think of a point in $\mathbb{R}^n$ as giving edge weight assignments for the edges of $G$. A point in $\mathbb{R}^n$ is in the tropical variety of the graphic matroid arising from $G$ if and only if each cycle in $G$ attains its minimum edge weight at least twice. We restrict our attention to graphic matroids arising from graphs with no loops and no parallel edges.

\begin{theorem} 

The unique minimal tropical basis for a graphic matroid is the collection of its \emph{induced cycles}.
\end{theorem}

An {\it induced cycle} is an induced subgraph that is itself a cycle.

\begin{proof}
Induced cycles are sufficient: Suppose we have a weight vector which achieves its minimum twice on each induced cycle. Consider a circuit arising from a non-induced cycle $C$. Then the induced subgraph on the vertices contained in $C$ contains a chord. Divide $C$ into cycles $C_1$ and $C_2$ along this chord. Continuing in this manner, we can decompose $C$ into a collection of induced cycles pasted together.  Thus by Lemma \ref{pasting}, induced cycles form a tropical basis.

Induced cycles are necessary: Suppose we have an induced cycle $C$ which is not in our tropical basis. Then it is possible to construct an edge weighting for which the minimum is attained at least twice on every cycle except $C$. Assign all but one edge of $C$ weight 1. Assign weight 0 to the remaining edge of $C$ and to all other edges of the graph. Every cycle incident with $C$ has two or more edges that are not contained in $C$, since otherwise $C$ would contain a chord and thus not be an induced cycle. Thus all cycles except $C$ have their minimum attained at least twice, while $C$ has a unique minimum. Thus it is necessary to have $C$ in every tropical basis.

\end{proof}

\subsection{Cographic Matroids}

A cographic matroid is formed from a connected graph $G$, having the edges of $G$ as its ground set. Circuits of cographic matroids are the inclusion-minimal edge cuts, i.e., sets of edges such that removing them makes $G$ disconnected.  In order to insure every circuit contains at least 3 elements, we will restrict our attention to cographic matroids arising from $3$-edge-connected graphs.  A graph is called \emph{$k$-edge-connected} if it remains connected after removing any $k-1$ edges.
A \emph{bridge} of a connected graph is an edge whose removal disconnects the graph.  A connected graph is $2$-edge-connected if and only if it does not contain a bridge.
 
 \begin{figure}[htb]
\includegraphics[scale=0.6]{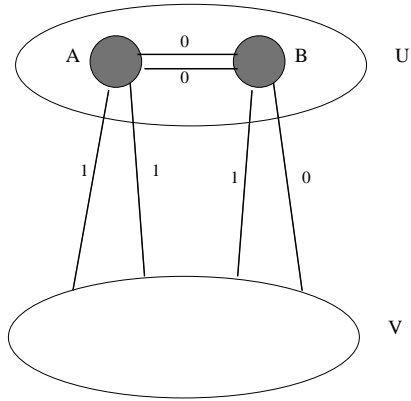}
\caption{2-edge-connected graph}
\label{2conngraph}
\end{figure}
  
\begin{theorem} The unique minimal tropical basis of a cographic matroid $M$ on graph $G$ consists of the edge cuts that split $G$ into two 2-edge-connected subgraphs.
\end{theorem}

\begin{proof}

For an edge cut $C$ that splits $G$ into connected subgraphs $U$ and $V$, let the \emph{index} of $C$ be the sum of the number of bridges in $U$ and $V$.  We will show by induction on the index that the circuits of positive index can be obtained by successively pasting index $0$ circuits, as in Lemma \ref{pasting}.

Suppose an edge cut $C$ splits $G$ into connected subgraphs $U$ and $V$.  Let $e$ be a bridge in $U$ that splits it into subgraphs $A$ and $B$.  Let $br_A, br_B$, and $br_V$ be the number of bridges in $A$, $B$, and $V$ respectively.  Then the index of $C$ is equal to $1 + br_A + br_B + br_V$.  Let $C'$ be the edge cut that splits $A$ and $B \cup V$.  Here $B \cup V$ denotes the induced subgraph of $G$ on the vertices in $B$ and $V$.  In $B \cup V$, there are at least two edges between $B$ and $V$ because otherwise the edge $e$ and the unique edge would disconnect $G$, contradicting the $3$-edge-connectedness of $G$.  Hence a bridge in $B \cup V$ must be either a bridge in $B$ or a bridge in $V$.  Therefore, the index of the edge cut $C'$ is at most $br_A + br_B + br_V$.  Similarly, the index of the edge cut $C''$ that splits $B$ and $A \cup V$ is at most $br_A + br_B + br_V$.   The edge cut $C$ is obtained by pasting the lower index circuits $C'$ and $C''$ along $e$.  By induction on the index, we see that $C$ is obtained by successively pasting index 0 circuits.  This proves (by Lemma \ref{pasting}) that the index 0 circuits form a tropical basis.

Now we show that all such edge cuts are necessary.
Suppose an edge cut $C \subset \text{edges}(G)$ splits the graph $G$ into subgraphs $U$ and $V$, each of them 2-edge-connected.  Consider edge weights as follows: each edge contained in $U$ and $V$ gets weight 0, one edge going between $U$ and $V$ gets weight 0, and the other edges between $U$ and $V$ get weight 1.  Consider any other edge cut $C' \neq C$.  It must cause either $U$ or $V$ to become disconnected.  Suppose it cuts $U$ into subgraphs $A$ and $B$ (see Figure \ref{2conngraph}). Since $U$ is 2-edge-connected, there are at least 2 edges going between $A$ and $B$, and those edges are in $C'$.  Hence the minimal weight $0$ is attained at least twice on $C'$, but not on $C$. Thus $C$ must be in our tropical basis.

\end{proof}

\subsection{The matroid $R_{10}$}

The matroid $R_{10}$ is a regular matroid whose elements are given by the edges of the complete graph on five vertices, $K_5$. Its circuits are given by the fifteen $4$-cycles of $K_5$, along with the complement of each $4$-cycle. Thus $R_{10}$ has 30 circuits.

\begin{proposition} The unique minimal tropical basis of $R_{10}$ consists of the fifteen 4-cycles.
\end{proposition}

\begin{proof}
The 4-cycles are necessary: For any 4-cycle $C=\{a_1,a_2,a_3,a_4\}$, we can construct a $0/1$-point $x$ which is excluded by $C$, but not any other circuit of $R_{10}$. Set $x_{a_1}=x_{a_2}=x_{a_3}=1$, and assign $0$ to all other entries of $x$. Then $x$ attains its minimum value on $C$ uniquely, while on every other circuit its minimum is attained at least twice. This is because no other $4$-cycle can share 3 elements with $C$, and any 4-cycle-complement will contain at least three elements disjoint from $C$, corresponding to `zero' entries of $x$.

In order to show the 4-cycles are sufficient, we show that any $0/1$-point excluded by a 4-cycle-complement will also be excluded by a 4-cycle. Consider a point $x$ which is excluded by a 4-cycle-complement $C=\{a_1,a_2,a_3,a_4,a_5,a_6\}$. Since $x$ attains its minimum uniquely on the terms of $C$, without loss of generality suppose $x_{a_1}=0$ and $x_{a_2}=x_{a_3}=x_{a_4}=x_{a_5}=x_{a_6}=1$. Note that no 4-cycle-complement contains a 4-cycle. In order for a 4-cycle to exclude $x$, it must intersect $C$ in exactly three elements other than $a_1$. All 4-cycle-complements look the same up to permuting the vertices. Figure \ref{r10} shows one such picture. One can easily see that for any 5-element subset of the 6 edges, there exists a 4-cycle meeting the subset in exactly 3 edges. Therefore every point excluded by a 4-cycle-complement is also excluded by a 4-cycle.

\begin{figure}[htb]
\includegraphics[scale=0.6]{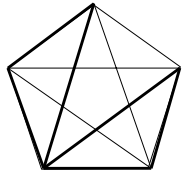}
\caption{4-cycle-complement}
\label{r10}
\end{figure}

\end{proof}

\begin{problem}
Do all regular matroids have unique minimal tropical bases?
\end{problem}

The Fano plane is not regular but has a unique minimal tropical basis.

\section{Non-constant coefficient case}

We now discuss the \emph{non-constant coefficient} case.  A {\it matroid polytope} of a rank $d$ matroid on the ground set $[n]$ is the convex hull of the characteristic vectors in $\RR^n$ of the bases of the matroid.  The {\it hypersimplex} $\Delta_{d,n}$ is the convex hull of $0/1$ vectors in $\RR^n$ with coordinate sum $d$.  A {\it matroidal subdivision} of $\Delta_{d,n}$  is a regular subdivision in which every face is a matroid polytope.   It was shown in \cite{Sp} that a point $p \in \TP^{[{n] \choose d}-1}$ is a tropical Pl\"{u}cker vector if and only if it induces a matroidal subdivision of  $\Delta_{n,d}$.  Moreover, the corresponding tropical linear space is the polyhedral subcomplex of the dual of the subdivision, consisting of the cells dual to matroid polytopes of loop-free matroids, i.e.\ matroids without any one-element circuits.

A coordinate of $p$ being $\infty$ is equivalent to deleting the corresponding vertex from the hypersimplex $\Delta_{d,n}$.  The case when all coordinates of $p$ are either $0$ or $\infty$ is the \emph{constant coefficient} case discussed in Sections 2 and 3.  In this case, the subdivision contains a single matroid polytope, and the circuits of the corresponding matroid are precisely the supports of the circuits defined in the Introduction.  The tropical linear space of a matroid is also called the \emph{Bergman fan}.

Since the hypersimplex $\Delta_{d,n}$ lies on the hyperplane of points whose coordinates sum to $d$, the duals of the subdivisions, hence the tropical linear spaces, have a lineality space containing $\RR(1,\dots, 1)$, as seen above.  Since each cell in a matroidal subdivision is a matroid polytope, the star of any face in the tropical linear space is isomorphic to a constant coefficient tropical linear space.

Let $K$ be the field of Puiseux series $\sum_{a \in I}c_a x^a$ where $I$ is a  locally finite subset of $\RR$ with a least element and $c_a \in \CC$.   Let the degree map $deg : K \rightarrow \TT  $ be the map that sends an element in $K \backslash \{0\}$ to its leading (lowest) exponent of $x$ and sends $0$ to $\infty$.  It induces a map from a projective space over $K$ to a tropical projective space.
  Let $L$ be a $d$-dimensional vector subspace in $K^n$, and let $P \in \PP_K^{{[n] \choose d} -1}$ be its Pl\"{u}cker vector, i.e.\ the the vector of maximal minors of a matrix whose row space is $L$.  Then  deg$(P)$ is a tropical Pl\"{u}cker vector.  If a tropical Pl\"{u}cker vector arises in this way, then the corresponding tropical linear space coincides with the image $deg(L) \subset \TT^n$ and is called \emph{realizable}.  If there is a representative Pl\"{u}cker vector containing only complex numbers, then we get the realizable constant coefficient case.  For a linear form $f = a_1 X_1 + \cdots + a_n X_n$ with $a_i \in K$, let its \emph{tropicalization} be the tropical linear form $(\text{deg}(a_1) \odot x_1) \oplus \cdots \oplus (\text{deg}(a_n) \odot x_n)$.
The circuits defined in the Introduction are precisely the tropicalizations linear forms with minimal support in the orthogonal complement of $L$.  We will refer to a set of linear forms over $K^n$ as a \emph{tropical basis} if their tropicalizations define the tropical linear space.

The \emph{tropical determinant} of an $r \times r$ square matrix $A = [a_{ij}]$ is defined to be
$$
\bigoplus_{\sigma \in S_r} (a_{1\sigma_1} \odot \cdots \odot a_{r \sigma_r}), \text{ where $S_r$ is the symmetric group of order $r!$.}
$$
The \emph{tropical rank} of a matrix is the largest integer $r$ such that there exists a submatrix which is tropically non-singular, i.e.\ the minimum in the tropical determinant is not unique.

 In the special case of realizable tropical linear spaces all of whose tropical Pl\"{u}cker coordinates are non-$\infty$, we conjecture a criterion for a set of linear forms to form a tropical basis, which generalizes Theorem 5.3 in \cite{RST}.

 \begin{theorem}
 \label{rank}
Let $L$ be a $n-k$ dimensional linear subspace in $K^n$ all of whose Pl\"{u}cker coordinates are non-zero.  Let $M\in K^{m \times n}, m \geq k,$ be a matrix whose rows are non-zero elements in the orthogonal complement of $L$.  If the rows of $M$ form a tropical basis for $L$, then any $k$ columns of deg$(M)$ have tropical rank $k$. 
 \end{theorem} 
 
 \begin{proof}
 Suppose that there is a $m \times k$ submatrix $A$ of deg$(M)$ with tropical rank less than $k$.  Then by \cite[Theorem 5.5]{DSS} the \emph{Kapranov rank} of $A$ is less than $k$, which means that there is an $m \times k $ matrix $A'$ over $K$ with rank less than $k$ such that $A = \text{deg}(A')$.  Let $v \in K^k$ be a non-zero vector in the kernel of $A'$.  Then deg$(v)$ is in the tropical prevariety in $\TT^k$ defined by the $m$ rows of $A$.  We can augment deg$(v)$ to a vector in the tropical prevariety of the rows of deg$(M)$ by putting $\infty$ in the other $n-k$ coordinates.  The support of this vector has size at most $k$, however, the points in the tropical linear space have support size at least $k+1$ because of the hypothesis that all the tropical Pl\"{u}cker coordinates are finite.  Hence the prevariety is not equal to the tropical linear space, so the rows of $M$ do not form a tropical basis.
 \end{proof}
 
 \begin{conjecture}
 The converse of Theorem \ref{rank} holds.
 \end{conjecture}
 
  \begin{example}
The proposition and the conjecture do not apply when some of the tropical Pl\"{u}cker coordinates are $\infty$.  Consider the $2$-dimensional linear subspace of $K^4$ which is the kernel of the matrix
 $$ M = 
\left( 
  \begin{array}{cccc}
    1 & 0 & 1 & 1 \\ 
    0 & 1 & 1 & 1 \\ 
  1 & -1 & 0 & 0\\
  \end{array}
\right), \text{ whose degree is }
\left(
  \begin{array}{cccc}
    0 & \infty & 0 & 0 \\ 
    \infty & 0 & 0 & 0 \\ 
    0 & 0 & \infty & \infty \\ 
  \end{array}
\right).
$$
The circuits are precisely these three rows, hence the $0/1$ points in the tropical variety are $(0,0,0,1), (0,0,1,0), (1,1,0,0)$.  The last two columns of deg$(M)$ have tropical rank only $1$, but the rows of $M$ form a tropical basis.  Notice that the corresponding tropical Pl\"{u}cker coordinate is $\infty$.  The first two rows of the matrix deg$(M)$ have the same tropical rank as the whole matrix in every subset of columns, but they do not form a tropical basis.  Hence we cannot determine whether a set of linear forms is a tropical basis just from the tropical ranks. \qed
 \end{example}
 
In the cases when Theorem \ref{rank} and its converse are applicable, we would get an algorithm for checking if a given matrix is a tropical basis.
 
\section{Parametrizations of Tropical Linear Spaces}
\label{sec:param}

So far we have been looking at tropical linear spaces as intersections of tropical hyperplanes.  In this section, we will look at them as images of tropical linear maps.

Let $A$ be an $n \times d$ matrix over $K$ whose image (column space) is $L$.  Let $deg(A)$ be the matrix whose entries are the degrees of entries in $A$.  This matrix defines a tropical linear map $deg(A) : \TT^d \rightarrow \TT^n$, $v \mapsto deg(A) \odot v$, where the tropical matrix multiplication $\odot$ is defined by replacing sums with minima and products with sums in the evaluation of the ordinary matrix product.

For any such $A$, we have
 \begin{equation}
 \label{eq2}
 deg(L) = deg(im(A)) \supseteq im(deg(A))
 \end{equation}
The containment holds because the columns of $deg(A)$ are in the tropical linear space $deg(L)$, and so are their tropical linear combinations, since tropical linear spaces are closed under taking tropical linear combinations.   Similar expressions hold in much more generality for tropical varieties, as shown in \cite{PS}.  We are interested in knowing when the equality $deg(L) = im(deg(A))$ is attained.

 A {\em cocircuit} of the linear space $L$ is an element in $L$ whose support is minimal with respect to inclusion.  They are circuits of the orthogonal complement $L^\perp$, i.e. the linear forms with minimal support that .  Two cocircuits with the same support must be constant multiples of each other, since otherwise one coordinate can be cancelled to get a vector with a smaller support.
 
 \begin{lemma}
Every nonzero $v \in L$ can be written as $v = v^1 + \dots + v^d$ for some cocircuits $v^i \in L$ such that $deg(v) = deg(v^1) \oplus \dots \oplus deg(v^d)$.
 \end{lemma}
 
 \begin{proof}
If $v$ is a cocircuit, then we are done.  

Suppose not.  Let $u \in L$ be a cocircuit with $supp(u) \subset supp(v)$.  Then for $c\in K$ with large enough degree, we have $deg(c u) \geq deg(v)$ coordinatewise.  Pick such a $c$ so that $c u$ and $v$ coincide in at least one coordinate, i.e.\ $supp(v - cu) \subsetneq supp(v)$.  Let $v^1 = c u$.  Since $deg(v^1) \geq deg(v)$, we have $deg(v) = deg(v^1) \oplus deg(v - v^1)$. 

If $v - v^1$ is not a cocircuit, then we can repeat the same argument on $v - v^1$, which has a strictly smaller support.  We will eventually end up with cocircuits with the desired properties.
 \end{proof}
 
 \begin{theorem}
 The equation $deg(L) = im(deg(A))$ holds if and only if every cocircuit in $L$ is represented in $A$. 
  \end{theorem}
 
 \begin{proof}
The ``if'' direction follows immediately from previous Lemma.
For the ``only if'' direction, suppose there is a cocircuit $c \in L$ whose support is not represented in $A$.  Then $deg(c) \in deg(L) \backslash im(deg(A))$ since any vector in $im(deg(A))$ with finite coordinates in $supp(c)$ also has finite coordinates outside $supp(c)$.
\end{proof}

\noindent This theorem for the constant coefficient case appeared in \cite[Proposition 7.5]{DSS}.

\section{Acknowledgments}
This project started during Debbie Yuster's Fall 2005 visit to the UC Berkeley Mathematics Department.  We are very grateful to Bernd Sturmfels for suggesting the problem and giving much guidance and encouragement along the way.  We thank Mike Develin for discussions and comments on an earlier draft of this paper, Dave Bayer for suggestions especially regarding Section 3, and Jenia Tevelev for discussions regarding Section 5.  Josephine Yu was supported by the NSF Graduate Research Fellowship and the UC Berkeley Graduate Opportunity Program Fellowship while working on this project.

 \end{document}